\documentclass[a4paper, 10pt, reqno]{article}

\usepackage{amsmath, amsthm, amsfonts, amssymb}

\newtheorem{theorem}{Theorem}
\newtheorem{lemma}[theorem]{Lemma}
\newtheorem*{Khinchin}{Khinchin's inequalities}

\DeclareMathOperator{\sgn}{sgn}

\newcommand{\R}{\mathbb{R}}

\newcommand{\ipr}[2]{\left\langle #1, #2 \right\rangle}
\newcommand{\bigipr}[2]{\bigl\langle #1, #2 \bigr\rangle}
\newcommand{\Bigipr}[2]{\Bigl\langle #1, #2 \Bigr\rangle}

\newcommand{\norm}[1]{\left\lVert#1\right\rVert}
\newcommand{\bignorm}[1]{\bigl\lVert#1\bigr\rVert}
\newcommand{\Bignorm}[1]{\Bigl\lVert#1\Bigr\rVert}
\newcommand{\biggnorm}[1]{\biggl\lVert#1\biggr\rVert}

\newcommand{\abs}[1]{\left\lvert#1\right\rvert}

\newcommand{\Bigabs}[1]{\Bigl\lvert#1\Bigr\rvert}

\begin{document}
\title{Vertex degrees of Steiner Minimal Trees in $\ell_p^d$ and other smooth 
Minkowski spaces}
\date{}
\author{K.\ J.\ Swanepoel \\ Department of Mathematics and Applied Mathematics 
\\ University of Pretoria \\ 0002 Pretoria,\, South Africa\\ 
e-mail: \texttt{konrad@math.up.ac.za}}
\maketitle
\begin{abstract}
We find upper bounds for the degrees of vertices and Steiner points in Steiner 
Minimal Trees in the $d$-dimensional Banach spaces $\ell_p^d$ independent of 
$d$.
This is in contrast to Minimal Spanning Trees, where the maximum degree of 
vertices grows exponentially in $d$ (Robins and Salowe, 1995).
Our upper bounds follow from characterizations of singularities of SMT's due 
to Lawlor and Morgan (1994), which we extend, and certain 
$\ell_p$-inequalities.
We derive a general upper bound of $d+1$ for the degree of vertices of an SMT 
in an arbitrary smooth $d$-dimensional Banach space (i.e.\ Minkowski space); 
the same upper bound for Steiner points having been found by Lawlor and Morgan.
We obtain a second upper bound for the degrees of vertices in terms of 
$1$-summing norms.
\end{abstract}
\section{Introduction}
Given a metric space $(X,\rho)$ and a set $S\subseteq X$, a {\em Minimal 
Spanning Tree (MST)} of $S$ is a tree $T$ with vertex set $V(T)=S$ and edge 
set $E(T)$ such that 
$$\sum_{\{x,y\}\in E(T)}\rho(x,y)$$
is minimal among all trees on $S$.

A {\em Steiner Minimal Tree (SMT)} of $S$ is a tree $T$ with vertex set $V(T)$ 
satisfying $S \subseteq V(T)\subseteq X$ such that
$$\sum_{\{x,y\}\in E(T)}\rho(x,y)$$
is minimal among all trees on $S$ with vertex sets satisfying 
$S \subseteq V(T)\subseteq X$.
The elements of $S$ are {\em vertices}, and the elements of $V(T)\setminus S$ 
are {\em Steiner points} of the SMT.

Estimates for the largest degrees of MST's and SMT's have consequences for the 
complexities of algorithms that find such trees. 
For example, it is known that an MST on $n$ points can be calculated in 
polynomial time \cite{CT}, while calculating the SMT in the euclidean or 
rectilinear planes is NP-hard \cite{GGJ, GJ}.
Upper bounds for the degrees of vertices and Steiner points are used to reduce 
the search space of known exponential time algorithms.

Distance functions other than euclidean or rectilinear are sometimes used.
The $\ell_p$ metrics have been found useful; see \cite{LoMo}.
We consider general Minkowski spaces, i.e.\ finite dimensional Banach spaces, 
and then specialize to $\ell_p^d$, $d$-dimensional real linear space with norm
$$\bignorm{(x_1,\dots,x_d)}_p = \biggl(\sum_{i=1}^d \abs{x_i}^p\biggr)^{1/p}.$$

It is known that in a Minkowski space, the largest degree of an MST is equal 
to the so-called Hadwiger number $H(B)$ of the unit ball $B$ of the space 
\cite{Cieslik}.
For each $1\leq p \leq\infty$ there is an exponential lower bound for the 
Hadwiger number of $\ell_p^d$, $H(B_p^d) > (1+\epsilon_p)^d$ \cite{RS}.

In contrast to this, we show in Section~\ref{results} that the degrees of 
both vertices and Steiner points of an SMT in $\ell_p^d$ ($1<p<\infty$) are 
bounded above by functions of $p$ alone, independent of $d$.
For $p>2$ we derive a general upper bound of $7$, with various sharper values 
for specific $p$.
For $1<p<2$ however, we find an upper bound exponential in $p^\ast:=p/(p-1)$, 
and a lower bound linear in $p^\ast$, as $p$ tends to $1$.
Thus with respect to the SMT problem, $\ell_p^d$ behaves very similarly to 
euclidean space, where both vertices and Steiner points have degree at most 3.

For general $d$-dimensional smooth Minkowski spaces, it is known that the 
degree of a Steiner point is at most $d+1$ \cite{LM}.
In Section~\ref{minkowski} we show that this upper bound also holds for the 
degree of a vertex in an SMT.
The proof has two ingredients.
Firstly, in Section~\ref{singularities} we derive a characterization of the 
local structure of a vertex in an SMT (Theorem~\ref{vertices}) similar to the 
characterization of Steiner points due to Lawlor and Morgan \cite{LM}.
We also rederive their characterization, paying attention to some 
combinatorial subtleties (Theorem~\ref{steinerpoints}).
Both derivations are completely elementary.
The second ingredient is Theorem~\ref{bound}, which generalizes a result of 
\cite{FLM} and \cite{LM}, thus answering a question in \cite{S}.

In Theorem~\ref{summing} we also obtain an upper bound for the degrees of 
vertices and Steiner points in terms of the $1$-summing norm of the dual of 
the space.

\section{Derivation of the singularity characterizations}\label{singularities}
Theorem~\ref{steinerpoints} below, due to \cite{LM}, provides a 
characterization of the structure of the neighbourhood of a Steiner point in 
an SMT in a smooth Minkowski space.
We give a similar characterization of the structure of the neighbourhood of a 
vertex in an SMT in Theorem~\ref{vertices}.
Both characterizations are in terms of unit vectors in the dual of the 
Minkowski space.

We now recall some facts about dual spaces.
Note that the discussion below pertains to finite dimensional Banach spaces, 
i.e.\ Minkowski spaces; see \cite{Thompson}.

For any $d$-dimensional real vector space $X$, the {\em dual} of $X$, denoted 
by $X^\ast$, is the vector space of linear functionals $x^\ast:X\to\R$.
This dual is also a $d$-dimensional vector space.
We denote application of $x^\ast\in X^\ast$ to $x\in X$ by $\ipr{x^\ast}{x}$.
If $X$ is furthermore a Minkowski space with norm $\norm{\cdot}$, then 
$\norm{x^\ast}^\ast = \sup_{\norm{x}\leq 1}\ipr{x^\ast}{x}$ defines a norm 
on $X^\ast$.

We say that a Minkowski space is {\em smooth} if 
$$\lim_{t\to 0} \frac{\norm{x+th}-\norm{x}}{t} =: f_x(h)$$
exists for all $x,h\in X$ with $x\neq 0$.
It follows easily that $f_x\in X^\ast$, $\norm{f_x}^\ast=1$ and 
$\ipr{f_x}{x}=\norm{x}$.
A linear functional $x^\ast\in X^\ast$ is a {\em norming functional} of $x$ 
if $x^\ast$ satisfies $\ipr{x^\ast}{x}=\norm{x}$ and $\norm{x^\ast}^\ast=1$.
Each non-zero vector in a Minkowski space has a norming functional (the 
Hahn-Banach theorem).
A Minkowski space is smooth iff each non-zero vector has a unique norming 
functional.

A Minkowski space $X$ is {\em strictly convex} if $\norm{x}=\norm{y}=1$ and 
$x\neq y$ imply that $\norm{\frac{1}{2}(x+y)}<1$, equivalently, that the 
boundary of the unit ball of $X$ does not contain any straight line segment.
A Minkowski space $X$ is smooth [strictly convex] iff $X^\ast$ is strictly 
convex [smooth].

The balancing and collapsing conditions in Theorems~\ref{steinerpoints} and 
\ref{vertices} thus occur in a strictly convex space.
We say that a finite set of unit vectors $x_1,\dots,x_m\in X$ satisfies the 
{\em balancing condition} if
\begin{equation}\label{balancing}
\sum_{i=1}^m x_i = 0,
\end{equation}
and satisfies the {\em collapsing condition} if
\begin{equation}\label{collapsing}
\Bignorm{\sum_{i\in J}x_i} \leq 1 \text{ for each } J\subseteq\{1,\dots,m\}.
\end{equation}
Note that the above balancing condition is the characterization of the 
so-called Fermat point of a set of points in a smooth Minkowski space in the 
non-absorbing case (i.e.\ where the Fermat point differs from the given 
points) in terms of norming functionals, derived in \cite{CG}.

\begin{theorem}[Lawlor and Morgan \cite{LM}]\label{steinerpoints}
Let $a_1,\dots,a_m$ be distinct non-zero points in a smooth Minkowski space 
$X$. For each $i=1,\dots, m$, let $a_i^\ast$ be the norming functional of 
$a_i$.
Then the tree connecting each $a_i$ to $0$ is an SMT of 
$S = \{a_1,\dots,a_m\}$ iff
$\{a_1^\ast,\dots,a_m^\ast\}$ satisfies the balancing and collapsing 
conditions in $X^\ast$.
\end{theorem}

\proof
$\Rightarrow$: Since we have an SMT, for any $x\in X$
$$\sum_{i=1}^m\norm{a_i-x}\geq\sum_{i=1}^m\norm{a_i},$$
i.e.\ for any unit vector $e\in X$ the function
$$\phi_e(t) := \sum_{i=1}^m (\norm{a_i+te}-\norm{a_i}) \geq 0$$
attains a minimum at $t=0$.
For sufficiently small $t$, $a_i+te\neq 0$, and $\phi_e(t)$ is differentiable 
at 0, with $\phi'_e(0)=0$.
But $$\phi'_e(0)=\lim_{t\to 0}\sum_{i=1}^m \frac{\norm{a_i+te}-\norm{a_i}}{t} 
= \sum_{i=1}^{m}\ipr{a_i^\ast}{e}.$$
Therefore, $\sum_{i=1}^ma_i^\ast = 0$.

Secondly, given $J\subseteq\{1,\dots,m\}$, define a tree $T_J$ as follows: 
Connect $\{a_i:i\in J\}$ to an arbitrary point $x$, connect 
$\{a_i:i\notin J\}$ to 0, and connect $x$ to $0$.
Then the total length of $T_J$ is not smaller than $\sum_{i=1}^m\norm{a_i}$:
$$\sum_{i\in J}\norm{a_i-x}+\sum_{i\notin J}\norm{a_i}+\norm{x} 
\geq\sum_{i=1}^m\norm{a_i},$$
i.e.\ for any unit vector $e$ the function 
$$\psi_e(t) := \sum_{i\in J}(\norm{a_i-te}-\norm{a_i})+\abs{t}\geq 0$$
attains a minimum at $t=0$.
However, $\psi_e$ is not differentiable at $0$.
Circumventing this difficulty, we calculate
\begin{align*}
0 & \leq \lim_{t\to 0^+}\frac{\psi_e(t)}{t} 
= \lim_{t\to 0^+}\sum_{i\in J}\frac{\norm{a_i-te}-\norm{a_i}}{t}+1 \\
  & = \sum_{i\in J}\ipr{a_i^\ast}{-e} + 1
\end{align*}
and $\ipr{\sum_{i\in J}a_i^\ast}{e} \leq 1$ for all unit $e$.
Thus $\norm{\sum_{i\in J}a_i^\ast}^\ast\leq 1$.

$\Leftarrow$: Let $a_1^\ast,\dots,a_m^\ast\in X^\ast$ satisfy 
\eqref{balancing} and \eqref{collapsing}, and let $T$ be any SMT of 
$\{a_1,\dots,a_m\}$.
We have to show that 
$$\sum_{\{x,y\}\in E(T)}\norm{x-y}\geq\sum_{i=1}^m \norm{a_i}.$$
For $i\geq 2$, let $P_i$ be any non-overlapping path in $T$ from $a_1$ to 
$a_i$, i.e.\ $P_i = x_{1}^{(i)}x_2^{(i)}\dots x_{k_i}^{(i)}$ with 
$x_1^{(i)}=a_1, x_{k_i}^{(i)}=a_i$ and $\{x_j^{(i)},x_{j+1}^{(i)}\}$ distinct 
edges in $E(T)$ for $j=1,\dots,k_i-1$.
Note that each edge of $T$ is used in some $P_i$, since the union of the 
paths is a connected subgraph of $T$.
For each edge $e\in E(T)$ we assign a direction depending on the way $e$ is 
traversed in some $P_i$ containing $e$.
This direction is unambigious, since if two paths would give conflicting 
directions, their union would contain a cycle.
We denote a directed edge from $x$ to $y$ by $(x,y)=\vec{e}$ and the set of 
directed edges by $\vec{E}(T)$.
For each $\vec{e}\in \vec{E}(T)$, let 
$S_{\vec{e}}:=\{i\geq 2:\vec{e}\in P_i\}$.
Then 
\begin{equation*}
\begin{split}
\sum_{i=1}^m \norm{a_i} & = \sum_{i=1}^m \ipr{a_i^\ast}{a_i} \\
& = \sum_{i=2}^m \ipr{a_i^\ast}{a_i-a_1} \qquad \text{ (by the balancing 
condition)} \\
& = \sum_{i=2}^m \sum_{j=2}^{k_i-1}\bigipr{a_i^\ast}{x_{j+1}^{(i)} 
-x_{j}^{(i)}} \\
& = \sum_{\vec{e}=(x,y)\in \vec{E}(T)} 
\sum_{i\in S_{\vec{e}}}\ipr{a_i^\ast}{y-x} \\
& \leq \sum_{\vec{e}=(x,y)\in \vec{E}(T)} 
\Bignorm{\sum_{i\in S_{\vec{e}}}a_i^\ast}^\ast\norm{x-y} \\
& \leq \sum_{(x,y)\in \vec{E}(T)} \norm{x-y} 
\qquad \text{ (by the collapsing condition)}.
\end{split}
\end{equation*}
\qed

As mentioned in \cite{LM}, the balancing and collapsing conditions are still 
sufficient for the tree in the above theorem to be an SMT in non-smooth 
spaces, if \eqref{balancing} and \eqref{collapsing} holds for {\em some} 
norming functional $a_i^\ast$ for each $a_i$.
A similar remark holds for the next theorem.

\begin{theorem}\label{vertices}
Given points $a_1,\dots,a_m\neq 0$ in a smooth Minkowski space $X$, let 
$a_i^\ast$ be the norming functional of $a_i$.
Then the tree connecting each $a_i$ to $0$ is an SMT of 
$S = \{0,a_1,\dots,a_m\}$ iff
$\{a_1^\ast,\dots,a_m^\ast\}$ satisfies the collapsing condition in $X^\ast$.
\end{theorem}

\proof
Similar to the proof of the previous theorem.
Note that there is no balancing condition, since we cannot perturb $0$, as 
$0$ is in this case a vertex of the SMT.
\qed

\section{Upper bounds for smooth Minkowski spaces}\label{minkowski}
For a Minkowski space $X$, let $v(X)$ be the largest degree of a vertex of an 
SMT in $X$, and $s(X)$ the largest degree of a Steiner point in an SMT.

In \cite{LM} it is shown that $s(X)\leq d+1$ if $X$ is smooth and 
$d$-dimensional.
This inequality is sharp in the sense that there are spaces and SMT's where 
the degree of $d+1$ is attained.
We give a similar bound for $v(X)$:

\begin{theorem}\label{smoothbound}
For a smooth Minkowski space $X$ of dimension $d\geq 2$, 
$$3\leq s(X)\leq v(X)\leq d+1.$$
The outer inequalities are sharp in general.
\end{theorem}

\proof
Theorems~\ref{steinerpoints} and \ref{vertices} immediately imply 
$s(X)\leq v(X)$.

In any $2$-dimensional subspace of the dual $X^\ast$ we can find two unit 
vectors $x^\ast,y^\ast$ such that $\norm{x^\ast-y^\ast}^\ast = 1$.
Then the set $\{x^\ast,-y^\ast, y^\ast-x^\ast\}$ satisfies \eqref{balancing} 
and \eqref{collapsing}.

The euclidean spaces $X=\ell_2^d$ are examples where $s(X)=v(X)=3$.

The rest of the theorem now follows from Theorem~\ref{vertices} and 
Theorem~\ref{bound} below.
An example where $v(X)=d+1$ may be constructed in the same way as for $s(X)$, 
as is done in \cite[Lemma 4.3]{LM}.
\qed

The following theorem, suggested in \cite{S}, sharpens results from 
\cite{FLM} and \cite{LM} by eliminating the balancing condition from the 
hypotheses.

\begin{theorem}\label{bound}
Let $X$ be a strictly convex $d$-dimensional Minkowski space.
If $x_1,\dots,x_m\in X$ are unit vectors satisfying the collapsing condition, 
then $m\leq d+1$.
Furthermore, if the balancing condition is not satisfied, i.e.\ 
$\sum_{i=1}^mx_i\neq 0$, then $m\leq d$.
\end{theorem}

\proof
Let $x_i^\ast\in X^\ast$ be norming functionals of $x_i$.
Firstly, for $i\neq j$ we have $$1+\ipr{x_i^\ast}{x_j} 
= \ipr{x_i^\ast}{x_i+x_j}\leq\norm{x_i+x_j}\leq 1$$
by the collapsing condition, and thus $$\ipr{x_i^\ast}{x_j}\leq 0 
\text{ for } i\neq j.$$

Secondly,
$$0\leq\Bigipr{x_i^\ast}{-\sum_{j\neq i} x_j}
\leq\biggnorm{\sum_{j\neq i} x_j}\leq 1.$$
If $\bigipr{x_i^\ast}{-\sum_{j\neq i} x_j} = 1$, then $x_i^\ast$ is also a 
norming functional of $-\sum_{j\neq i} x_j$, which is now a unit vector.
Then, since $X$ is strictly convex, it easily follows that 
$x_i=-\sum_{j\neq i} x_j$.

Thus, if $\sum_{i=1}^m x_i\neq 0$, then
$$0\leq\Bigipr{x_i^\ast}{-\sum_{j\neq i} x_j}<1,$$
and the diagonal of the matrix $A=\bigl[\ipr{x_i^\ast}{x_j}\bigr]_{i,j=1}^m$ 
majorizes the rows.
Thus $A$ is invertible.
Since $A$ has rank at most $d$, we obtain $m\leq d$.

If however $\sum_{i=1}^m x_i = 0$, the above argument applied to 
$x_1,\dots,x_{m-1}$ gives $m-1\leq d$.
\qed

Note that in the above proof, we do not nearly use the full force of the 
collapsing condition.

For the next bound, we recall a notion from the local theory of Banach spaces.
The {\em absolutely summing constant} or the {\em $1$-summing norm} (of the 
identity operator on) a Minkowski space $X$ is defined to be
$$\pi_1(X):= \inf\Bigl\{c>0:\forall x_1,\dots,x_m\in X:\;\sum_{i=1}^m
\norm{x_i}\leq c\max_{\epsilon_i=\pm 1}
\Bignorm{\sum_{i=1}^m\epsilon_ix_i}\Bigr\}.$$
This notion has been studied extensively; see e.g.\ 
\cite{Mc, DR, R, KS, GG, KT}.
Note that the quantity $(2\pi_1(X))^{-1}$ has also been called the 
{\em Macphail constant} in the literature.

\begin{theorem}\label{summing}
For a smooth Minkowski space $X$,
$$s(X)\leq v(X)\leq 2\pi_1(X^\ast).$$
\end{theorem}

\proof
Let $x_1^\ast,\dots,x_m^\ast\in X^\ast$ be unit vectors satisfying the 
collapsing condition, with $m=v(X)$.
Then, for any sequence of signs $\epsilon_i=\pm 1,i=1,\dots,m$ we have 
$\norm{\sum_i\epsilon_i x_i^\ast}^\ast\leq 2$, hence
$$m=\sum_{i=1}^m\norm{x_i^\ast}^\ast 
\geq \frac{m}{2}\max_{\epsilon_i=\pm 1}
\Bignorm{\sum_{i=1}^m\epsilon_i x_i^\ast}^\ast,$$
implying that $\frac{m}{2} \leq \pi_1(X^\ast)$.
\qed

It is known that $\sqrt{d}\leq\pi_1(X)\leq d$ for any $d$-dimensional $X$ 
\cite{KS}.
We thus obtain an upper bound worse than that of Theorem~\ref{smoothbound}, 
although it is of the same order.
It is however possible in principle to obtain bounds better than that of 
Theorem~\ref{smoothbound} for specific spaces.
However, we cannot do better than $2\sqrt{d}$.

\section{Upper bounds for $\ell_p^d$}\label{results}
Restricting ourselves to the smooth case $1<p<\infty$, we recall that the 
dual of $\ell_p^d$ is $(\ell_p^d)^\ast = \ell_{p^\ast}^d$, where 
$1/p + 1/p^\ast = 1$.
We use the Khinchin inequalities with the best constants, due to 
\cite{Szarek} and \cite{Haagerup1,Haagerup2}.

\begin{Khinchin}
For any $1\leq q < \infty$ there exist constants $A_q, B_q > 0$ such that for 
any $a_1,\dots,a_n\in\R$ we have
$$A_q\biggl(\sum_{i=1}^na_n^2\biggr)^{1/2} 
\leq \biggl(2^{-n}\sum_{\epsilon_{i}=\pm 1}
\Bigabs{\sum_{i=1}^{n}\epsilon_i a_i}^q\biggr)^{1/q} 
\leq B_q\biggl(\sum_{i=1}^na_n^2\biggr)^{1/2}.$$
For $q\geq 2$ we have $A_q=1$, $B_q=\sqrt{2}\bigl(\Gamma(\frac{q+1}{2})/
\sqrt{\pi}\bigr)^{1/q}$, and for $1\leq q \leq 2$, $B_q=1$, 
$$A_q=\begin{cases}2^{1/2-1/q} \text{ if } q < q_0, \\
\sqrt{2}\bigl(\Gamma(\frac{q+1}{2})/\sqrt{\pi}\bigr)^{1/q} \text{ if } q 
\geq q_0,
\end{cases}$$
where $q_0 \approx 1.8474$ is defined by $\Gamma\left(\frac{q_0+1}{2}\right)
=\frac{\sqrt{\pi}}{2}$, $1<q_0<2$.
\qed
\end{Khinchin}

The following lemma is analogous to \cite[Hilfsatz 4]{CS}.
We omit the proof, which easily follows from calculus.

\begin{lemma}\label{pineq}
Let $x,y\in\R$ and $1\leq q\leq 2$.
Then $$\abs{x+y}^q\geq 2^{q-2}\bigl(\abs{x}^{q/2}\sgn{x} 
+ \abs{y}^{q/2}\sgn{y}\bigr)^2.$$
\end{lemma}
\qed

The earliest reference we could find to the following lemma is Rankin 
\cite{Rankin}.

\begin{lemma}\label{inner}
Let $x_1,\dots,x_m\in\ell_2^d$ satisfy $\norm{x_i}_2=1$ and 
$\ipr{x_i}{x_j} < -1/n$ for $i\neq j$, where $n$ is a positive integer.
Then $m\leq n$.
\end{lemma}

\proof
\begin{equation*}
\begin{split}
0 & \leq\Bignorm{\sum_{i=1}^mx_i}_2^2 = \sum_{i=1}^m\norm{x_i}_2^2 
+ 2\sum_{i<j}\ipr{x_i}{x_{j}} \\
  & < m - m(m-1)/n.
\end{split}
\end{equation*}
\qed

The next two theorems show that the largest degree of a vertex $v(\ell_p^d)$ 
and the largest degree of a Steiner point $s(\ell_p^d)$ in an SMT in 
$\ell_p^d$ are both relatively small and independent of $d$.
In particular, for $p\geq 2$ we have a general upper bound of $7$.
For $2\leq p \lesssim 3.40942$ we furthermore obtain the exact values of 
$v(\ell_p^d)$ and $s(\ell_p^d)$.
For $p<2$ we only obtain a lower bound linear in $p^\ast$ and an upper bound 
exponential in $p^\ast$.
It is not clear what the correct order of growth should be in this case.

\begin{theorem}
Let $2\leq p<\infty$ and $d\geq 3$.
\begin{align}
&s(\ell_p^d) = v(\ell_p^d) = 3 \text{ for } 2\leq p 
< \tfrac{\log 3}{\log 3 - \log 2} \approx 2.70951, \label{one}\\
&s(\ell_p^d) = v(\ell_p^d) = 4 \text{ for } \tfrac{\log 3}{\log 3 - \log 2} 
\leq p < \tfrac{\log 8 - \log 3}{\log 4-\log 3}\approx 3.40942,\label{two}\\
&4\leq s(\ell_p^d)\leq v(\ell_p^d) \leq 5 \text{ for } 
\tfrac{\log 3}{\log 3 - \log 2} \leq p 
< \tfrac{\log 4}{\log 4-\log 3}\approx 4.81884,\label{three}\\
&4\leq s(\ell_p^d)\leq v(\ell_p^d) \leq 6 \text{ for } 
\tfrac{\log 3}{\log 3 - \log 2} \leq p 
< \tfrac{\log 4}{\log 8-\log 7}\approx 10.3818,\label{four}\\
&4\leq s(\ell_p^d)\leq v(\ell_p^d) \leq 7 \text{ for all } p 
\geq \tfrac{\log 3}{\log 3 - \log 2}.\label{five}
\end{align}
\end{theorem}

\proof
Let $q:= p^\ast = p/(p-1)$.
The lower bound of $3$ for $s(X)$ and $v(X)$ comes from 
Theorem~\ref{smoothbound}.
For $p\geq (\log 3)/(\log 3 - \log 2)$, i.e.\ for $q\leq (\log 3)/(\log 2)$, 
we obtain $4$ unit vectors in $\ell_q^d$ satisfying the balancing and 
collapsing conditions as follows:
\begin{align*}
& x_1:= 3^{-1/q}(1,1,1), & x_2:= 3^{-1/q}(1,-1,-1), \\
& x_3:= 3^{-1/q}(-1,1,-1), & x_4:= 3^{-1/q}(-1,-1,1).
\end{align*}
 
For the upper bounds, let $x_1,\dots,x_m\in\ell_q^d$ be unit vectors 
satisfying the collapsing condition.

We first use a ``twisting'' technique used in the Geometry of Numbers; see 
\cite{Rineq}.
Denote the coordinates of $x_i$ as $x_i = (x_{i,1},x_{i,2},\dots,x_{i,d})$.
Define $\tilde{x}_i=(\tilde{x}_{i,1},\tilde{x}_{i,2},\dots,\tilde{x}_{i,d})$ 
by $\tilde{x}_{i,n}:=\abs{x_{i,n}}^{q/2}\sgn{x_{i,n}}$.
Note that $\norm{\tilde{x}_i}_2=1$, i.e.\ we have twisted $x_i$ to become a 
euclidean unit vector.
By Lemma~\ref{pineq} we obtain for $i\neq j$ that
$$1\geq\norm{x_i+x_j}_q^q\geq 2^{q-2}\norm{\tilde{x}_i+\tilde{x}_j}_2^2 
= 2^{q-2}(2+2\ipr{\tilde{x}_i}{\tilde{x}_j}),$$
where $\ipr{\cdot}{\cdot}$ denotes the standard euclidean inner product.
Thus $\ipr{\tilde{x}_i}{\tilde{x}_j}\leq 2^{1-q}-1 < 0$.
If $p<(\log 3)/(\log 3-\log 2)$, i.e.\ $q>(\log 3)/(\log 2)$, then 
$2^{1-q}-1<-\frac{1}{3}$.
By Lemma~\ref{inner} we obtain $m\leq 3$, and \eqref{one} follows.
Similarly, if $p<(\log 8 - \log 3)/(\log 4-\log 3)$, then 
$2^{1-q}-1<-\frac{1}{4}$, hence $m\leq 4$, and \eqref{two} follows.

For the remaining estimates we apply Khinchin's inequalities.
We may assume in the light of \eqref{one} and \eqref{two} that 
$p\geq (\log 8 - \log 3)/(\log 4-\log 3)$, i.e.\ 
$q \leq  (\log 8-\log 3)/(\log 2) < q_0$.
Thus $A_q = 2^{1/2-1/q}$.
By \eqref{collapsing} we have for any sequence of signs 
$\epsilon_i=\pm 1, i=1,\dots,m$ that 
$\norm{\sum_{i=1}^m\epsilon_ix_i}_q\leq 2$.
Therefore,
\begin{align*}
2^q & \geq \sum_{n=1}^d 2^{-m}\sum_{\epsilon_i=\pm 1}
\Bigabs{\sum_{i=1}^m\epsilon_{i}x_{i,n}}^q \\
    & \geq \sum_{n=1}^d A_q^q\biggl(\sum_{i=1}^m x_{i,n}^2\biggr)^{q/2}
&\text{ (Khinchin's inequality)} \\
    & = A_q^q\sum_{n=1}^d\bignorm{\bigl(\abs{x_{i,n}}^q\bigr)_i}_{2/q}
&\text{ (where $\bigl(\abs{x_{i,n}}^q\bigr)_{i=1}^m\in\ell_{2/q}^m$)}\\
    & \geq A_q^q\biggnorm{\Bigl(\sum_{n=1}^d\abs{x_{i,n}}^q\Bigr)_i}_{2/q}
&\text{ (triangle inequality in $\ell_{2/q}^m$)} \\
    & = A_q^q\biggl(\sum_{i=1}^m\norm{x_i}_q^2\biggr)^{q/2} = A_q^q m^{q/2},
\end{align*}

and $m\leq 4/A_q^2 = 2^{3-2/p}<8$.
Estimates \eqref{three}, \eqref{four} and \eqref{five} now follow.
\qed

\begin{theorem}
Let $1<p<2$ and $d\geq 3$.
Then
\begin{equation}
\min(d,f(p^\ast))\leq s(\ell_p^d), v(\ell_p^d) \leq \min(d+1,2^{p^\ast}),
\label{six}
\end{equation}
where for $q>2$,
$$f(q):= \max\{d: 2(d-2)^q+(d-2)2^q\leq(d-1)^q+d-1\}.$$
In particular, 
\begin{align*}
& f(q)\geq 3 \text{ for } q> 2, \\
& f(q)\geq 4 \text{ for } q\geq 3.21067, \\
& f(q)\geq 5 \text{ for } q\geq 3.40093, \\
& f(q)\geq \lceil q/\log 2\rceil \text{ for } q\geq 3.69247.
\end{align*}
\end{theorem}

\proof
Let $q:= p^\ast = p/(p-1)$.

The upper bound follows from Theorem~\ref{smoothbound} and an application of 
Khinchin's inequalities:
\begin{align*}
2^q & \geq \sum_{n=1}^d 2^{-m}\sum_{\epsilon_i=\pm 1}
\Bigabs{\sum_{i=1}^m\epsilon_i x_{i,n}}^q \\
    & \geq \sum_{n=1}^d \biggl(\sum_{i=1}^m x_{i,n}^2\biggr)^{q/2}
&\text{ (Khinchin's inequality)} \\
    & = \sum_{n=1}^d\norm{x_i}_2^q \\
    & \geq \sum_{n=1}^d\norm{x_{i}}_q^q = m
&\text{ (monotonicity of $q$-norms)}.
\end{align*}

For the lower bound we may assume that $d\geq 4$.
Let $x_i$ be the vector in $\ell_q^d$ with $d-1$ in its $i$'th coordinate, 
and $-1$ in the remaining coordinates, for $i=1,\dots,d$.
Let $\hat{x}_i:= \norm{x_i}_q^{-1}x_i$.
Then $\{\hat{x}_i:i=1,\dots,d\}$ satisfies the balancing condition 
\eqref{balancing}.
This set will also satisfy the collapsing condition iff for all 
$2\leq k\leq d/2$,
\begin{equation*}
g(k,d,q) := k(d-k)^q + (d-k)k^q \leq (d-1)^q+d-1 = g(1,d,q).
\end{equation*}
By differentiating with respect to $q$ and using $2\leq k\leq d/2$, it is 
easily seen that if $g(k,d,q)\leq g(1,d,q)$ holds for some $q=q'$, then it 
will hold for all $q\geq q'$.
The following numerical facts are easily verified:
\begin{align*}
& g(k,d,q)\leq g(2,d,q)\text{ for }4\leq d\leq 7,\, 2\leq k\leq d/2 
\text{ and } p\geq 3.2,  \\
& g(2,4,q)\leq g(1,4,q) \text{ for }q\geq 3.21066\dots, \\
& g(2,5,q)\leq g(1,5,q) \text{ for }q\geq 3.40092\dots, \\
& g(2,6,q)\leq g(1,6,q) \text{ for }q\geq 3.69246\dots, \text{ and }\\
& g(2,7,q)\leq g(1,7,q) \text{ for }q\geq 4.09345\dots.
\end{align*}
It is now sufficient to show for $d\geq 8$ and $q=(d-1)\log 2$ that 
$g(k,d,q)\leq g(2,d,q)\leq g(1,d,q)$ for all $2\leq k \leq d/2$. 
Firstly, note that in this case $g(2,d,q)\leq g(1,d,q)$ is equivalent to
$$2^{1+(d-1)\log(d-2)} + (d-2)2^{(d-1)\log 2} \leq 2^{(d-1)\log(d-1)} + d-1,$$
which is easily verified for $d\geq 8$.

Secondly, to show that $g(k,d,q)\leq g(2,d,q)$ it is sufficient to show that 
$$f(x):=x(1-x)^q+(1-x)x^q, \qquad\tfrac{2}{d} \leq x \leq \tfrac{1}{2}$$
attains its maximum at $x=\frac{2}{d}$.
To see this, it is in turn sufficient to show that $f'(x)\leq 0$ for 
$2/d\leq x\leq 1/2$.
By setting $y=(1-x)/x$ we find that it is sufficient to show that for 
$1\leq y\leq d/2-1$,
$$x^{-q}f'(x) = y^q-qy^{q-1}-1+qy =: h(y) \leq 0.$$
By calculating the first and second derivatives of $h(y)$ and recalling that 
$q>3$, it is seen that $h(y)$ does not attain its maximum if $1<y<d/2-1$.
Since $h(1)=0$, we only have to show that $h(d/2-1)\leq 0$, which easily 
follows from $q\geq 4$ and $d\geq 8$.
\qed

\section*{Acknowledgement}
This paper is part of the author's PhD thesis being written under supervision 
of Prof.\ W. L. Fouch\'e at the University of Pretoria.

\end{document}